\documentclass[12pt]{article}
\usepackage{amssymb}
\def\b#1{{\bf #1}}
\def\c#1{{\cal #1}}

\def\1{{\bf 1}}
\newcommand{\botimes}{\mbox{\boldmath \boldmath$\otimes$}}
\newcommand{\tl}{\,\triangleleft}
\def\cross{{\triangleright\!\!\!<}}
\def\cocross{{>\!\!\!\triangleleft\,}}
\def\g{\mbox{\bf g\,}}

\def \uqg{\mbox{$U_q{\/\mbox{\bf g}}$ }}
\def\R{{\cal R}\,}

\newtheorem{theorem}{Theorem}

\newcommand{\be}{\begin{equation}}
\newcommand{\ee}{\end{equation}}
\newcommand{\bea}{\begin{eqnarray}}
\newcommand{\eea}{\end{eqnarray}}
\newcommand{\ba}{\begin{array}}
\newcommand{\ea}{\end{array}}

\begin{document}

\title{Decoupling of Tensor factors in Cross Product and Braided Tensor
Product Algebras\footnote{Contribution to the Proceedings of the
``International 
Colloquium on Group Theoretical Methods in Physics'' (Group24), 
Paris, July 2002.}}
\author{Gaetano Fiore}

\author{        Gaetano Fiore, \\\\
         \and
        Dip. di Matematica e Applicazioni, Fac.  di Ingegneria\\ 
        Universit\`a di Napoli, V. Claudio 21, 80125 Napoli
        \and
        I.N.F.N., Sezione di Napoli,\\
        Complesso MSA, V. Cintia, 80126 Napoli
        }
\maketitle
\date{}

\begin{abstract}
We briefly review and illustrate
 our procedure to `decouple' by transformation of generators:
 either a Hopf algebra $H$ 
from a $H$-module algebra $\c{A}_1$ in their cross-product
$\c{A}_1\cocross H$; or
two (or more) $H$-module algebras $\c{A}_1,\c{A}_2$. 
These transformations are based on the existence 
of an algebra map $\c{A}_1\cocross H\to \c{A}_1$.
\end{abstract}

\vfill
\noindent
Preprint 02-64 Dip. Matematica e Applicazioni, Universit\`a di Napoli\\
DSF/29-2002
\newpage

\section{Decoupling of tensor factors in
cross product algebras}

Let $H$  be a Hopf algebra over $\b{C}$, say,  $\c{A}$ a unital
(right, say) $H$-module algebra. We denote by $\tl$  the
right action, namely the bilinear map
such that, for any $a,a'\in\c{A}$ and
$g,g'\in H$
\be
\ba{c}
\tl: (a,g)\in\c{A} \times H\to a\tl g\in\c{A},  \\
a\tl(gg') = (a\tl g) \tl   g',           \qquad  \qquad
(aa')\tl g = (a\tl g_{(1)})\, (a'\tl g_{(2)}).
\ea       \label{modalg}
\ee
We have used the Sweedler-type notation
 $\Delta(g)=g_{(1)}\botimes g_{(2)}$ 
for the coproduct $\Delta$.
The cross-product algebra $\c{A} \cocross H$ is $H\otimes\c{A}$
as a vector space,
and so we denote as usual $g\otimes a$ simply by $ga$;
$H\1_{\c{A}}$, $\1_H\c{A}$ are subalgebras
isomorphich to $H, \c{A}$, and so we omit to 
write either unit $\1_{\c{A}},\1_H$
whenever multiplied by non-unit elements;  
for any $a\in\c{A}$, $g\in H$ the product fulfills
\be
a g=g_{(1)}\, (a\tl g_{(2)}).                \label{crossprod}
\ee
$\c{A} \cocross H$ is a $H$-module algebra  itself, if we
extend $\tl$ on $H$ as the adjoint action:
$h\tl g= Sg_{(1)}\, h\, g_{(2)}$.
If $H$ is a Hopf $*$-algebra, $\c{A}$ a $H$-module
$*$-algebra,  then,  as known,
these two $*$-structures can be glued in a unique one to make
$\c{A}\cocross H$ a $*$-algebra  itself.

\begin{theorem}\cite{Fio01} 
Let $H$ be a Hopf algebra, $\c{A}$ a right
$H$-module algebra. If there exists a "realization" $\tilde\varphi$
of $\c{A} \cocross H$ within $\c{A}$ acting as
the identity on $\c{A}$, i.e. an algebra map
\be
\:\:\tilde\varphi: \c{A}\cocross H \rightarrow \c{A},
\qquad\qquad\qquad\qquad\tilde\varphi(a)=a,               \label{ident0r}
\ee
then $\tilde\zeta(g):=g_{(1)} \tilde\varphi(S g_{(2)})$ defines
an injective algebra map $\tilde\zeta: H\to \c{A}\cocross H$ such that 
\be
[\tilde\zeta(g), \c{A}]=0
\ee
for any $g\in H$. Moreover $\c{A}\cocross H=\tilde\zeta(H)\,\c{A}$.
Consequently, 
the center of the cross-product algebra
$\c{A}\cocross H$ is given by
$\c{Z}(\c{A}\cocross H)=\c{Z}(\c{A})\,\tilde\zeta\left(\c{Z}(H)\right)$,
and if $H_c,\c{A}_c$ are Cartan subalgebras of $H$
and $\c{A}$ respectively, then $\c{A}_c\,\tilde\zeta(H_c)$
is a Cartan subalgebra of $\c{A}\cocross H$.
Finally, if $\tilde\varphi: \c{A}\cocross H \rightarrow \c{A}$
is a $*$-algebra map, then also
$\tilde\zeta: H \rightarrow \tilde\c{C}$ is.
\end{theorem}
\label{theorem1r}

The equality $\c{A}\cocross H=\tilde\zeta(H)\,\c{A}$ means that
$\c{A}\cocross H$ is equal to a product $\c{A}\, H'$,
where $H'\equiv\tilde\zeta(H)\subset \c{A}\cocross H$ is a subalgebra
isomorphic to $H$ and {\it commuting}
with $\c{A}$, i.e. $\c{A}\cocross H$ 
is isomorphic to the {\it ordinary} tensor product
algebra $H \otimes \c{A}$. In other words, if 
$\{a_I\},\{g_J\}$ are resp. sets of generators of $\c{A}, H$, 
then $\{a_I\}\cup\{\tilde\zeta(g_J)\}$ is a more manageable
set of generators of $\c{A}\cocross H$ than $\{a_I\}\cup\{g_J\}$.
The other statements allow to determine
Casimirs and complete sets of commuting observables, key ingredients
to develop representation theory.

Well-known examples of maps (\ref{ident0r}) are the ``vector field'' or the
``Jordan-Schwinger'' realizations of the UEA $U\g$'s
($\g$ being a Lie algebra),
where $\c{A}$ is resp.
the Heisenberg algebra on a $\g$-covariant space or  a $\g$-covariant
Clifford algebra. In the case of e.g. the
Heisenberg algebra on the Euclidean space $\b{R}^3$, the well-known
realization of the three generators $J^{ij}$, ($i\neq j$) of $\g=so(3)$
as ``vector fields'' (namely homogeneous first order differential operators)
$$
\tilde\varphi(J^{jk}):=x^jp^k-x^kp^j
$$
gives nothing but the orbital angular momentum operator
in 1-particle  quantum mechanics
($x_i,p^i\equiv-i\partial^i$ denote the position and momentum
components respectively). Then
$$
\tilde\zeta(J^{jk})=J^{jk}-\tilde\varphi(J^{jk})= J^{jk}-(x^jp^k-x^kp^j).
$$
gives the difference between the total and the orbital angular
momentum, i.e. the ``intrinsic'' angular momentum (or ``spin''), which
indeed commutes with the $x^i,p^i$'s.
Maps $\tilde\varphi$ have been determined \cite{ChuZum95,Fiocmp95} 
also for a number of
\uqg-covariant, i.e. quantum group covariant,
deformed Heisenberg (or Clifford) algebras.
So the theorem is immediately
applicable to them.

No map $\tilde\varphi$ can exist if  we take as $\c{A}$ just the space
on which the Heisenberg algebra is built (in the previous example
$\b{R}^3$), since
one cannot realize the non abelian algebra $\c{A}\cocross H$ in terms
of the abelian $\c{A}$.
Surprisingly, a map $\tilde\varphi$ may exist if we deform the algebras.
For instance, in Ref. \cite{GroMadSte00} a map $\tilde\varphi$
realizing  $U_qso(3)$ has been determined
for $\c{A}$ the $q$-deformed fuzzy sphere $\hat S^2_{q,M}$.
To treat other examples
we generalize the previous results by
weakening our assumptions. Namely, we require at least
that $H$ admits a Gauss decomposition
$$
H=H^+H^-=H^-H^+
$$
into two Hopf subalgebras $H^+,H^-$
for each of which analogous maps $\tilde\varphi^+,\tilde\varphi^-$
(\ref{ident0r}) coinciding on $H^+\cap H^-$ exist.
Then Theorem \ref{theorem1r} will apply separately to
$\c{A}\cocross H^+$ and $\c{A}\cocross H^-$.  What about
the whole $H\cross \c{A}$?

\begin{theorem}\cite{Fio01}
Under the above assumptions setting
$\tilde\zeta^{\pm}(g^{\pm}):=g^{\pm}_{(1)}\tilde\varphi^{\pm}(S  
g^{\pm}_{(2)})$
(where $g^{\pm}\in H^{\pm}$ respectively)
defines injective algebra maps 
$\tilde\zeta^{\pm}: H^{\pm}\to \c{A}\cocross H^{\pm}$ such that  
$[\tilde\zeta^{\pm}(g^{\pm}), \c{A}]=0$ for any $g^{\pm}\in H^{\pm}$.
Moreover
\be
\c{A} \cocross H=\tilde\zeta^+(H^+)\,\tilde\zeta^-(H^-)\,\c{A}
=\tilde\zeta^-(H^-)\,\tilde\zeta^+(H^+)\,\c{A}.         \label{decom}
\ee
Any $c\in\c{Z}(\c{A}\cocross H)$ can be expressed in the form
\be
c=\tilde\zeta^+\left(c^{(1)}\right)\tilde\zeta^-\left(c^{(2)}\right)c^{(3)},
\label{centre+-}
\ee
where $c^{(1)}\otimes c^{(2)}\otimes c^{(3)}
\in H^+\otimes H^- \otimes\c{Z}(\c{A})$ and
$c^{(1)}c^{(2)}\otimes c^{(3)}
\in \c{Z}(H) \otimes\c{Z}(\c{A})$; conversely any
such $c$ $\in\c{Z}(\c{A}\cocross H)$.
If $H_c\subset H^+ \cap H^-$ and
$\c{A}_c$ are Cartan subalgebras resp. of $H$
and $\c{A}$, then $\c{A}_c\,\tilde\zeta^+(H_c)$
[$\equiv\c{A}_c\,\tilde\zeta^-(H_c)$]
is a Cartan  subalgebra of $\c{A}\cocross H$.

If $\tilde\varphi^{\pm}$ are $*$-algebra map, then also
$\tilde\zeta^{\pm}: H^{\pm} \rightarrow \tilde\c{C}$ are.
If
$\tilde\varphi^{\pm}\Big((\alpha^{\mp})^*\Big)=
[\tilde\varphi^{\mp}(\alpha^{\mp})]^*$
 $\forall \alpha^{\mp}\in \c{A} \cocross H^{\mp} $,
then $\tilde\zeta^{\pm}(g^{\mp}{}^*)=[\tilde\zeta^{\mp}(g^{\mp})]^*$,
with $g^{\mp}\in H^{\mp}$.
\label{theorem+-r}
\end{theorem}

 As a consequence of this theorem,  $\forall g^+\in H^+$, $g^-\in H^-$
$\exists c^{(1)}\otimes c^{(2)}\otimes c^{(3)}\in
\c{Z}(\c{A})\otimes H^-\otimes H^+$
 (depending on $g^+,g^-$) such that
\be
\tilde\zeta^+(g^+)\tilde\zeta^-(g^-)=c^{(1)}\tilde
\zeta^-(c^{(2)})\tilde\zeta^+(c^{(3)}).            \label{com+-}
\ee
 These will be the "commutation relations"  between elements of
$\tilde\zeta^+(H^+)$ and $\tilde\zeta^-(H^-)$. Their form will depend on
the specific algebras considered.

As an application we consider now the pair $(H,\c{A})$
with $H=U_qso(3)$ and $\c{A}=\b{R}_q^3$, the
(algebra of functions on) the 3-dim quantum Euclidean space,
whose generators we denote resp. by $E^+,E^-,K,K^{-1}$ and $p^+,p^0,p^-$.
In our present conventions
\bea
&& p^0p^{\pm}=q^{\pm 1}p^{\pm}p^0\qquad\qquad [p^+,p^-]=(1-q^{-1})\,p^0p^0
\label{Arel}\\
&& K\,E^{\pm}=q^{\pm 1}E^{\pm}K\qquad\qquad [E^+,E^-]_{q^{-1}}=\frac{K^2-1}
{q^2-1}\label{Hrel}\\
&& \ba{lll}
Kp^0=p^0K,\qquad & Kp^{\pm}=q^{\mp 1}p^{\pm}K\qquad
& [p^0,E^{\pm}]=\mp p^{\mp}\\
p^{\pm}E^{\mp}=q^{\pm 1}E^{\mp}p^{\pm}\qquad & [p^+,E^+]_q=p^0\qquad
& [p^-,E^-]_{q^{-1}}=\!-\!q^{-1}\!p^0
\ea \label{AHrel}\\
&& \Delta(K)=K\otimes K\qquad\qquad
\Delta(E^{\pm})=E^{\pm}\otimes K+\1\otimes E^{\pm} \label{coprod}\\
&& (p^0)^*=p^0, \quad (p^-)^*=p^+, \quad K^*=K,\quad  (E^+)^*=E^-,
\label{star3}
\eea
where $[a,b]_w:=ab-wba$. The first three relations
give the algebra structure of
$\b{R}_q^3\cocross U_qso(3)$ (this underlies the quantum
group of inhomogenous transformations of $\b{R}_q^3$), 
(\ref{coprod}) together with
$\varepsilon(E^{\pm})=0$, $\varepsilon(K)=1$ the coalgebra of $H$,
(\ref{star3}) the $*$-structure corresponding to compact $H$ and ``real''
$\b{R}_q^3$ (this requires $q\in\b{R}$).
The element $P^2:=qp^+p^-+p^0p^0+p^-p^+$ is central, 
positive definite under this
$*$-structure and real under the other one (that requires $|q|=1$).
We enlarge the algebra by introducing also the square root and the
inverse $P,P^{-1}$. Setting $P=1$ we obtain the quantum Euclidean
sphere $S_q^2$, and $S_q^2\cocross U_qso(3)$ can be interpreted as
the algebra of observables of a quantum particle on $S_q^2$.
The maps $\tilde\varphi^+,\tilde\varphi^-$,
$\tilde\zeta^+,\tilde\zeta^-$, the algebra relations
among the new generators
$e^+:=\tilde\zeta^+(E^+)$, $e^-:=\tilde\zeta^-(E^-)$,
$k:=\tilde\zeta^+(K)=\tilde\zeta^-(K)$, and the additional central
element $c$ of the form (\ref{centre+-}) are given by
\bea
&& \tilde\varphi^{\pm}(K)=\eta\frac P{p_0}\qquad
\tilde\varphi^+(E^+)=\frac 1{(q\!-\!1)p^0}p^-\qquad
\tilde\varphi^-(E^-)=\frac q{(q\!-\!1)p^0}p^+\\
&& k=K\frac{\eta p_0}{P}\qquad
e^+=\frac \eta{ P}
\left[E^+p^0 +\frac {qp^-}{1\!-\!q}\right]\qquad
e^-=\frac \eta{ P}\left[E^-p^0 +\frac {p^+}{1\!-\!q}\right]\\
&& k\, e^{\pm}=q^{\pm 1}
e^{\pm} k\qquad
[e^+,e^-]_{q^{-1}}=
\frac{k^2+ \eta^2}{q^2-1}  \qquad c=\frac{e^+e^-k^{-1}}{1\!+\!q}\!+\!
\frac{k\!-\!qk^{-1}}
{(1\!-\!q^2)^2}. \label{Zrel}
\eea
Here $\eta\in\b{C}$, $\eta\neq 0$. (\ref{Zrel}) translates (\ref{com+-}),
and differs from (\ref{Hrel})$_3$ by the presence at the rhs of the
``central charge'' $(\eta^2+1)/(q^2-1)$. Only for $\eta^2=-1$
can the maps $\tilde\varphi^+,\tilde\varphi^-$ and
$\tilde\zeta^+,\tilde\zeta^-$ be glued into maps $\tilde\varphi$
and $\tilde\zeta$ respectively; then the latter will be
$*$-maps under the non-compact $*$-structure where $|q|=1$, but not
under (\ref{star3}). In order this to happen we need to take
$\eta\in\b{R}$. The $*$-representations of the $e^+,e^-,k$ subalgebra
for real $q$ differ from the ones of $U_qsu(2)$ in that they
are lowest-weight but not highest-weight representations, or viceversa
\cite{Fio95}.

In Ref. \cite{FioSteWes00,Fio01} we give maps $\tilde\varphi^{\pm}$
for the cross products $\b{R}_q^N\cocross U_qso(N)$ for all $N\ge 3$.

\section{Decoupling of braided tensor products}

As known,
if $H$ is a noncocommutative Hopf algebra (e.g. a quantum group $\uqg$)
and $\c{A}_1, \c{A}_2$ are two (unital)
$H$-module algebras the tensor product algebra
$\c{A}_1\botimes \c{A}_2$ will not
be in general a $H$-module algebra. If $H$ is quasitriangular
a $H$-module algebra can be obtained as the
{\it braided} tensor product algebra
$\c{A}^+:=\c{A}_1\underline{\otimes}^+\c{A}_2$, which
is defined as follows:
the vector space
underlying the latter is still the tensor product
of the vector spaces underlying $\c{A}_1, \c{A}_2$,
and so we shall denote as usual $a_1\otimes a_2$ simply by $a_1a_2$;
$\c{A}_1\1_{\c{A}_2}$, $\1_{\c{A}_1}\c{A}_2$
are still subalgebras
isomorphich to $\c{A}_1, \c{A}_2$, and so one can omit to 
write the units, whenever they are mutiplied by non-unit
elements; but the ``commutation relation'' between
$a_1\in\c{A}_1, a_2\in\c{A}_2$ are modified:
\be
a_2a_1= (a_1\tl \R^{(1)})\, (a_2\tl \R^{(2)}).     \label{absbraiding}
\ee
Here $\R\equiv \R^{(1)}\otimes \R^{(2)}\in H^+\otimes H^-$
(again a summation symbol at the rhs has been suppressed)
denotes the so-called universal $R$-matrix or
quasitriangular structure of $H\equiv$
\cite{Dri86},
and as before $H^{\pm}$ denote some positive and negative
Borel Hopf subalgebras of $H$.
If $H=\uqg$, then in the limit $q\to 1$ $H$ becomes
the cocommutative Hopf algebra $U\g$ and
$\R\to \1\otimes\1$. As a consequence $a_2a_1\to a_1a_2$
and thus $\c{A}^+$ goes to the ordinary tensor product algebra.
An alternative braided tensor product
$\c{A}^-=\c{A}_1\underline{\otimes}^-\c{A}_2$.
can be obtained by replacing in (\ref{absbraiding})
$\R$ by $\R^{-1}_{21}$. This is equivalent to exchanging
$\c{A}_1$ with $\c{A}_2$.

$\c{A}^+$ (as well as $\c{A}^-$) is a $*$-algebra
if $H$ is a Hopf $*$-algebra, $\c{A}_1$, $\c{A}_2$ are $H$-module
$*$-algebras (we use the same symbol $*$ for the $*$-structure
on all algebras $H,\c{A}_1$, etc.), and 
$\R^*\equiv\R^{(1)}{}^*\otimes \R^{(2)}{}^* =\R^{-1}$.
In the quantum group case this requires $|q|=1$.
Under the same assumptions also $\c{A}_1\cocross H$ is
a $*$-algebra.

If $\c{A}_1, \c{A}_2$ represent the algebras of
observables of different two quantum systems, (\ref{absbraiding})
will mean that in the composite system their degrees of freedom
are "coupled" to each other.
But again one  can ``decouple'' them by a transformation of generators if there
exists an algebra map $\tilde\varphi_1^+$, or
an algebra map $\tilde\varphi_1^-$:

\begin{theorem} \cite{FioSteWes00}.
Let $\{H,\R\}$ be a quasitriangular Hopf algebra and
$H^+,H^-$ be Hopf subalgebras of $H$ such that $\R\in H^+\otimes H^-$.
Let $\c{A}_1, \c{A}_2$ be respectively a $H^+$- and
a $H^-$-module algebra, so that we can define $\c{A}^+$
as in (\ref{absbraiding}), and $\tilde\varphi_1^+$ be a map
of the type (\ref{ident0r}), so that we can define
the ``unbraiding'' map $\chi^+:\c{A}_2\rightarrow \c{A}^+$ by
\begin{equation}
\chi^+(a_2):= \tilde\varphi_1^+(\R^{(1)})\, (a_2\tl \R^{(2)}).
\label{Def+}
\end{equation}
Then $\chi^+$ is an injective algebra map such that
\begin{equation}
[\chi^+(a_2),\c{A}_1]=0,                       \label{commute}
\end{equation}
namely the subalgebra
$\tilde\c{A}_2^+:=\chi^+(\c{A}_2)\approx\c{A}_2$ commutes
with $\c{A}_1$. Moreover $\c{A}^+=\c{A}_1\tilde\c{A}_2^+$.
Finally, if $\R^*=\R^{-1}$ and $\tilde\varphi_1^+$
is a $*$-algebra map then $\chi^+$  is, and 
$\c{A}_1$, $\tilde\c{A}_2^+$ are closed under $*$.

By replacing everywhere $\R$ by $\R^{-1}_{21}$ we obtain
an analogous statement valid for $\c{A}^-,\chi^-$.
\label{maintheo}
\end{theorem}

Of course,
we can use the above theorem  iteratively to completely
decoupe the  braided tensor product algebra of an arbitrary number $M$ of
copies of $\c{A}_1$. One can also combine the two methods illustrated
here to decouple the tensor factors in `mixed' tensor products such as
\be
(\c{A}_1\cocross H)\,\underline{\otimes}^{\pm} \c{A}_2,\qquad
\c{A}_1\underline{\otimes}^{\pm} (\c{A}_2\cocross H), \qquad
(\c{A}_1\underline{\otimes}^{\pm} \c{A}_2)\cocross H.\qquad\qquad 
                                                    \label{combined}
\ee

\end{document}